\documentclass[11pt]{amsart}
\usepackage[english]{babel}
\usepackage{amsfonts, amsmath, amsthm, amssymb,amscd,indentfirst}

\numberwithin{equation}{section}

\usepackage{marginnote}

\usepackage{mathrsfs}
\usepackage{euscript}

\usepackage{xcolor}
\usepackage{amsmath,amssymb,latexsym,indentfirst}
\usepackage{times}
\usepackage{palatino}
\usepackage{stmaryrd}

\usepackage{wasysym}

\usepackage{graphicx}
\graphicspath{ {./images/}}

\theoremstyle{plain}
\newtheorem{theorem}{Theorem}[section]

\newtheorem{proposition}{Proposition}[section]

\theoremstyle{definition}
\newtheorem{definition}{Definition}[section]

\newtheorem{example}{Example}[section]

\theoremstyle{remark}
\newtheorem{remark}{Remark}[section]

\newcommand{\interior}[1]{%
	{\kern0pt#1}^{\mathrm{o}}%
}

\def\dirac{{\partial\!\!\!/}}

\def\Tc{{\mathcal T}}
\def\eh{{{\bf{\rm{e}}}}}

\def\Rc{{\mathcal R}}

\def\Fc{{\mathcal F}}
\def\Ec{{\mathcal E}}
\def\Uc{{\mathcal U}}

\def\Hc{{\mathcal H}}
\def\Dc{{\mathcal D}}
\def\Fc{{\mathcal F}}
\def\Gc{{\mathcal G}}
\def\Hc{{\mathcal H}}

\def\bg{{\mathcal{B}\mathcal{G}}}

\def\ka{{K_{{{\rm area}}}}}

\title{Infinite connected sums, 
	$K$-area and positive scalar curvature}%

\author{Levi Lopes de Lima}
\address{Universidade Federal do  Cear\'a,
	Departamento de Matem\'atica, Campus do Pici, R. Humberto Monte, s/n, 60455-760,
	Fortaleza/CE, Brazil}
\thanks{Partially supported by CNPq/Brazil grant
	312485/2018-2 and by FUNCAP/CNPq/PRONEX grant 00068.01.00/15.}
\email{levi@mat.ufc.br}

\begin{document}
	\maketitle

	\begin{abstract}
		Whyte \cite{whyte2001index} used the index theory of Dirac operators
		and Block-Weib\-erger uniformly finite homology \cite{block1992aperiodic} to show that
		certain infinite connected sums do not carry a metric with
		nonnegative scalar curvature in their bounded geometry class. His
		proof uses a coarse version of the $\hat{A}$-class to obstruct such
		metrics. In this note we prove a version of Whyte's result where a variant of the notion of
		infinite $K$-area, originally  due to Gromov \cite{gromov1996positive}, is used to
		obstruct metrics with positive scalar curvature.
	\end{abstract}
	
	\section{Introduction}
	\label{intro}
	
	We  consider the category $\bg_n$ of manifolds with bounded
	geometry. Thus, objects in $\bg_n$ are complete $n$-dimensional
	Riemannian manifolds whose curvature tensor and covariant
	derivatives of all orders are uniformly bounded and whose
	injectivity radius is positive. The morphisms of $\bg_n$ are
	diffeomorphisms with bounded distortion. This allows us to split
	the space of such metrics into classes, the so-called {\em bounded
		geometry classes} of metrics. Whenever we refer to a manifold, it
	should be understood that it is equipped with a complete metric varying
	within a fixed bounded geometry class. Also, all manifolds in this
	note will be spin, unless otherwise stated. Finally, in what follows we assume that $n\geq 4$.
	
	We observe that in \cite{block1992aperiodic} a homology theory $H_0^{\rm uf}$, named
	uniformly finite homology in degree zero, has been defined which is
	preserved by the above morphisms and therefore is a bounded geometry
	invariant. In particular, a subset $S\subset X$ defines a homology
	class $[S]\in H^{\rm uf}_0(X)$ if it is {\it locally uniformly finite}
	in the sense that for each $r>0$ there exists $C_r>0$ such that the
	amount of points of $S$ inside any metric ball of radius $r$ is
	bounded from above by $C_r$ (for more on the functor $H^{\rm uf}_0$ and
	its relation to bounded de Rham cohomology, see Section \ref{bd}).
	
	Now take $Y\in\bg_n$, $M$ a closed manifold and $S\subset Y$ as above.
	Let $Y\sharp_{{S}}M$ be the manifold obtained by connected summing
	to $Y$ a copy of $M$ along a small neighborhood of each element of
	$S$.
	We remark that 
		$\bg_n$ is stable under such infinite connected sums, so
		that $Y\sharp_{{S}}M\in\bg_n$ as well.

	Whyte \cite{whyte2001index} used the Atiyah-Patodi-Singer index theory of Dirac operators to
	show that if $Y$ carries a metric of nonnegative scalar curvature
	(in a given bounded geometry class), $\widehat{A}(M)\neq 0$ and $[S]\neq
	0$ then $Y\sharp_{{S}}M$ does not carry a metric of nonnegative
	scalar curvature (in the corresponding bounded geometry class).
	
	\begin{example}\label{exam1}
		Let $S=\mathbb{Z}^{n}\subset\mathbb{R}^{n}$ be the standard integer
		lattice in flat Euclidean space. Then  $[\mathbb{Z}^{n}]\neq 0$ and
		hence $\mathbb{R}^{4l}\sharp_{\mathbb{Z}^{4l}}M^{4l}$ does not carry
		a metric with nonnegative scalar curvature if $\widehat{A}(M)\neq 0$
		(for example we can take $M$ to be the product of $l$ $K3$
		surfaces). More generally, we can take $Y\to Y_0$ to be an infinite
		covering with $Y_0$ closed and carrying a metric of nonnegative
		scalar curvature ($Y$ is equipped with the covering metric).  Then
		Whyte's result applies to $Y\sharp_{{S}}M$ if $S\subset Y$ is an
		orbit under deck transformations and $\pi_1(Y_0)$ is amenable (for
		more on Whyte's result, see Remarks \ref{handle} and \ref{new}).
	\end{example}
	
	The purpose of this note is to prove a version of Whyte's result
	that uses a notion of $K$-area which is well suited  
	to obstructing metrics with positive
	scalar curvature on these infinite connected sums. As explained in Section \ref{bd} below, this invariant, denoted $K_{\rm area}^b(X,g)$, is indeed a Riemannian invariant of $(X,g)\in \bg_{2k}$. However,  the fact that it vanishes, is finite or infinite is an invariant of the bounded geometry class of the metric. This allows us to simply write
	$K_{\rm area}^b(Y)=K_{\rm area}^b(Y,g)$, with the {\em proviso} that $g$ is allowed to vary in its bounded geometry class. 
	
	For our purposes, it is convenient to decompose $\bg_n=\bg_n^{A}\sqcup\bg_n^{N}$, where 
	\[
	\bg_n^{A}=\left\{X\in\bg_n;H_0^{\rm uf}(X)\neq \{0\}\right\},
	\]
	and $\bg_n^{N}=\bg_n\backslash \bg_n^{A}$. Elements in $\bg_n^{A}$ (respectively, $\bg_n^{N}$) are called {\em amenable} (respectively, {\em non-amenable}); see Proposition \ref{bg:known:0} for a geometric characterization of this splitting of $\bg_n$. With this terminology, our obstruction to positive scalar curvature goes as follows.
	
	\begin{theorem}\label{main}
	Let $Y\in\bg_{2k}^A$ be spin with $0\leq K_{\rm area}^b(Y)<+\infty$ and pick $H_0^{\rm uf}(Y)\ni [S]\neq 0$. Then
	$Y\sharp_{S}M$ carries no metric with positive scalar curvature if $M$ is spin with $K_{\rm area}(M)=+\infty$.
	\end{theorem}
	
	Here, $K_{\rm area}(M)$ is defined as in \cite[Section 4]{gromov1996positive}; see also Section \ref{kar} below for a review of this classical notion.

	Theorem \ref{main} follows immediately from the next two results.

	\begin{theorem}\label{aux2}
	Let $Y\in\bg_{2k}^A$ with $0\leq K_{\rm area}^b(Y)<+\infty$ and pick $H_0^{\rm uf}(Y)\ni [S]\neq 0$. Then
	$K_{\rm area}^b(Y\sharp_{S}M)=+\infty$ if $K_{\rm area}(M)=+\infty$.
\end{theorem}
	
	\begin{theorem}\label{aux1}
		If $X^{2k}\in \bg_{2k}$ is spin  with 	$K_{\rm area}^b(X)=+\infty$ then it
	carries no metric (in the given bounded geometry class) whose scalar
	curvature has non-positive part concentrated outside  arbitrarily
	large compact sets. In particular, it carries no metric with positive scalar curvature.
	\end{theorem}

	Theorem \ref{aux1} has an independent
		interest as it shows that infinite $K^b$-area is an obstruction to the
		existence of metrics of {positive} scalar curvature in the bounded
		geometry setting. In fact, it provides a somewhat stronger result, namely, given an exhaustion $W_i\subset X$ by compact domains then there is no sequence of metrics $g_i$ in the given bounded geometry class with the non-positive part of the scalar curvature of $g_i$ contained in $X\backslash W_i$.

	\begin{remark}\label{rem2}
		The assumption $\widehat{A}(M)\neq 0$ only makes sense if
		$\dim M$ (and hence $\dim Y$) is a multiple of 4. On the other
		hand, our result applies to  certain manifolds in every even
		dimension and moreover the attached manifold  $M$ can be chosen to
		be more familiar. It applies for example to $\mathbb R^{2k}
		\sharp_{\mathbb Z^{2k}}\mathbb{T}^{2k}$, where $\mathbb R^{2k}\in\bg_{2k}^A$ is the flat euclidean space, $\mathbb Z^{2k}\subset\mathbb R^{2k}$ is the standard integer lattice and  
		$\mathbb{T}^{2k}$ is a torus (which has infinite $K$-area). More generally, we can 		replace $\mathbb{T}^{2k}$ by any finitely enlargeable spin manifold
		\cite{gromov1980spin}; see Remark \ref{hyper}. Also, as observed in \cite{whyte2001index}, the class $[S]\in H_0^{\rm uf}(X)$, $X\in\bg_n^A$, 
		lies in a non-Hausdorff homology group and hence standard
		obstructions based on $C^*$-algebra techniques do not seem to work
		here.
	\end{remark}
	
	\begin{remark}\label{handle}
		It is shown in \cite{block1992aperiodic} that if $Y^{4l}$, $l\geq 2$, admits a
		metric with uniformly positive scalar curvature, $[S]= 0$ and $M$ is
		{\em simply connected} with $\hat{A}(M)\neq 0$ then $Y\sharp_S M$
		admits a positive scalar curvature metric. In view of this, Whyte's
		result means that $[S]\neq 0$ is the {\em only} obstruction to the
		existence of such metrics on $Y\sharp_S M$. The proof  in
		\cite{block1992aperiodic} uses a relative version of a surgery argument due to
		Gromov and Lawson \cite{gromov1980classification}, which by its turn is based on Smale's
		$h$-cobordism theory. In particular, simply connectedness of $M$
		seems to be essential there (in order to kill handles with small index)
		and consequently the argument breaks down in our case because
		$\ka(M)=+\infty$ implies the non-triviality of  $\pi_1(M)$ (see
		Remark \ref{hyper}).
	\end{remark}
	
	Our proof of Theorem \ref{main} follows Whyte's approach with
	suitable modifications to account for the  fact that we will be
	dealing with almost flat complex bundles over $X$. The argument
	makes  use of the index theory of
	Atiyah-Patodi-Singer (APS) type boundary conditions (see Section
	\ref{s2}). The classical
	notion of $K$-area is reviewed in Section \ref{kar} and our variant, which is well suited to the $\bg_{2k}$ category, is introduced in Section
	\ref{bd}. The proofs of Theorems \ref{aux2} and
	\ref{aux1} are presented in Section \ref{proofaux1}.

	\section{APS index theory} \label{s2}
	
	Let  $W$ be an oriented $n$-dimensional spin manifold with a fixed
	spin structure \cite{lawson2016spin}. In the presence of a Riemannian metric $g$, there exists
	over $W$ a canonical hermitian vector bundle $S_W$, the {spinor
		bundle}, which comes equipped with a Clifford product
	$\mathfrak c:\Gamma(TW) \to\Gamma({{\rm End}}(S_W))$ and a compatible
	connection $\nabla:\Gamma(S_W)\to \Gamma(T^*W\otimes S_W)$. Using
	these structures we can define the corresponding Dirac operator
	$\dirac:\Gamma(S_W)\to \Gamma(S_W)$ acting on spinors,
	\[
		\dirac=\sum_{i=1}^n \mathfrak c(\eh_i)\nabla_{\eh_i},
	\]
	where $\{\eh_i\}$ is a local orthonormal basis tangent to $W$. More
	generally, we can fix a hermitian vector bundle $\Ec$ with
	compatible connection $\nabla$ and consider the twisted Dirac
	operator $\dirac_{\Ec}:\Gamma(S_W\otimes\Ec)\to
	\Gamma(S_W\otimes\Ec)$ acting on (twisted) spinors.  The
	Weitzenb\"ock decomposition for the corresponding Dirac Laplacian is
	\begin{equation}\label{dec}
		\dirac_{\Ec}^2=
		\nabla^*\nabla+\frac{1}{4}\kappa+\Rc^{\left[\nabla\right]},
	\end{equation}
	where $\nabla^*\nabla$ is the Bochner Laplacian of $S_W\otimes\Ec$,
	$\kappa$ is the scalar curvature of $W$ and for
	$\psi\otimes\eta\in\Gamma(S_W\otimes\Ec)$,
	\begin{equation}\label{jos3}
		\Rc^{\left[\nabla\right]}(\psi\otimes\eta)=\frac{1}{2}\sum_{ij}\mathfrak c(\eh_i)\mathfrak c(\eh_j)\psi\otimes
		R^{\nabla}_{\eh_i,\eh_j}\eta,
	\end{equation}
	with $R^{\nabla}$ being the curvature tensor of $\nabla$. 
	
	If $W$ is
	closed, which we assume for the moment, $\dirac_{\Ec}$ is a self-adjoint elliptic operator and $\ker
	\dirac_{\Ec}$,  the space of  harmonic spinors, has finite
	dimension.
	If $n=2k$ one has a decomposition 
	\begin{equation}\label{decom}
		S_W\otimes \Ec=\left(S_W^+\otimes\Ec\right)\oplus
		\left(S_W^-\otimes\Ec\right)
	\end{equation}
	into positive and
	negative spinors induced by the chirality operator and
	$\dirac_{\Ec}$ interchanges the factors, so we can decompose,
	according to (\ref{decom}),
	\[
		\dirac_{\Ec}=\left(
		\begin{array}{cc}
			0 & \dirac^-_{\Ec} \\
			\dirac^+_{\Ec} & 0
		\end{array}
		\right),
	\]
	where
	\[
		\dirac^{\pm}_{\Ec}={\dirac_{\Ec}}|_{\Gamma(S_W^{\pm}
			\otimes\Ec)}:\Gamma(S_W^{\pm}\otimes\Ec)\to
		\Gamma(S_W^{\mp}\otimes\Ec).
	\]
	It follows that   $\dirac_{\Ec}^+$ and $\dirac_{\Ec}^-$
	are adjoint to each other, so we  get a well-defined index
	\[
		\label{index3} \hbox{ind}\, \dirac^+_{\Ec}=\dim\ker \dirac^+_{\Ec}
		-\dim\ker \dirac^-_{\Ec}.
	\]
	The Atiyah-Singer index formula computes this integer  as
	\begin{equation}\label{sigfor}
		{{\rm{ind}}}\, \dirac^+_{\Ec}=\int_W
		[\widehat{A}(TW)\wedge{{\rm{ch}}}(\Ec)]_{2k},
	\end{equation}
	where $\widehat{A}(TW)\in H^{4*}(W;\mathbb Q)$ is the $\widehat{A}$-class of
	$TW$, ${{\rm{ch}}}(\Ec)\in H^{2*}(W;\mathbb Q)$ is the Chern character of
	$\Ec$ and the notation $[\,\,\,\,]_{2k}$ means that integration
	picks the element of degree $2k$ in the wedge product. Specializing
	to the case where $k=2l$ and $\Ec=W\times\mathbb C$ is the trivial line
	bundle (equipped with a flat connection) we get
	
	\begin{equation}\label{sigfor2}
		{{\rm{ind}}}\, \dirac^+=\widehat{A}(W),
	\end{equation}
	where
	\begin{equation}\label{sigfor3}
		\widehat{A}(W)=\int_W [\hat{A}(TW)]_{4l}
	\end{equation}
	is the $\hat{A}$-genus of $W$. Notice that in this case (\ref{dec})
	reduces to
	\begin{equation}\label{reduc}
		\dirac^2=\nabla^*\nabla+\frac{\kappa}{4}.
	\end{equation}
	
	\begin{remark}\label{new}
		The  famous Lichnerowicz's argument \cite{lawson2016spin} is based
		on the fact that, since $\nabla^*\nabla$ is nonnegative, the
		positivity of $\kappa$ in (\ref{reduc}) implies that $\dirac$ is
		positive and hence invertible, which gives $\widehat{A}(W)=0$ by
		(\ref{index3}) and (\ref{sigfor2}). So, $\widehat{A}(W)\neq 0$ is a
		topological obstruction to the existence of metrics with $\kappa>0$.
		Thus, the point of Whyte's theorem is that if an obstructing $M$
		(i.e. with $\widehat{A}(M)\neq 0$) is `'glued'' to $Y$ along a 
		nontrivial class $[S]\in H^{\rm uf}_0(Y)$ then an obstruction to metrics with
		$\kappa\geq 0$ on $Y\sharp_{S}M$ arises, even if $Y$ originally
		carries such a metric. 
		Our main result (Theorem \ref{main}) says that
		instead of nonzero $\hat{A}$-{genus} we can use infinite $K$-area as
		a `'glued'' obstruction to metrics of {\em positive} scalar
		curvature as long as $0\leq K^b_{\rm area}(Y)<+\infty$.
	\end{remark}

	We now consider the index theory for manifolds with boundary (see
	for example \cite{atiyah1975spectral,booss1993bavnbek,gilkey1993index}).
	Assume from now on that $W$ is a compact spin manifold with
	dimension $n=2k$ and {\it non\-emp\-ty} smooth boundary
	$\Sigma\subset W$, and $\Ec$ is a hermitian bundle over $W$ with a
	compatible connection. Introduce
	Fermi coordinates $(x,u)\in \Sigma\times [0,\delta)\to \Uc$ in a
	collar neighborhood $\Uc$ of $\Sigma$ and set $\Sigma_u=\{(x,u)
	;x\in\Sigma\}$ so that $\Sigma_0=\Sigma$. Then, restricted to $\Uc$,
	\[
		\dirac=\mathfrak c\left({\partial_u}\right)\left({\partial_u}+ D-\frac{1}{2}H\right),
	\]
	where $H$ is the mean curvature of the embeddings
	$\Sigma_u\subset\Uc$ {(computed with respect to the inward unit vector field)} and  $D$ is a self-adjoint linear operator, the
	tangential Dirac operator, defined as follows. For each $u$,
	$S_W|_{\Sigma_u}$ comes equipped with the Clifford product
	$\mathfrak c^u=-\mathfrak c(\partial/\partial u)\mathfrak c$, so if we consider the
	induced connection
	$$
	\nabla^u=\nabla-\frac{1}{2}\mathfrak c^u(A),
	$$
	where $A$ is the shape operator of the embedding $\Sigma_u\subset
	\Uc$, then
	$$
	D=\sum_{i=1}^{2k-1}\mathfrak c^u(\eh_i)\nabla^u_{\eh_i},
	$$
	where $\{\eh_i\}$ is an orthonormal basis tangent to  $\Sigma_u$.
	After twisting with $\Ec$, we obtain a first order self-adjoint
	elliptic operator $D_{\Ec}$ acting on
	$\Gamma(S_{\Sigma}\otimes\Ec|_{\Sigma})$ and commuting with
	chirality, so we can decompose
	$S_W\otimes\Ec|_{\Sigma}=:\mathbb{S}_{\Ec}=\mathbb{S}^+_{\Ec}\oplus\mathbb{S}^-_{\Ec}$
	and accordingly,
	\[
		D_{\Ec}=\left(
		\begin{array}{cc}
			D_{\Ec}^+ & 0 \\
			0 & D_{\Ec}^-
		\end{array}
		\right)
	\]
	with $D_{\Ec}^{\pm}$ self-adjoint. Under the natural identification
	$\mathbb{S}^+_{\Ec}=\mathbb{S}^-_{\Ec}$ one has
	$D^+_{\Ec}=-D^-_{\Ec}$ and hence ${{\rm Spec}}(D_{\Ec})$ is
	symmetric with respect to $0\in\mathbb{R}$, but of course this does
	not need happen to the factors $D^{\pm}_{\Ec}$. Thus, for ${{\rm
			Re}}\,z\gg 0$ we define the eta function
	$$
	\eta_{\Ec}^+(z)=\sum_{0\neq\lambda\in {{\rm Spec}}(D_{\Ec}^+)} {{\rm
			sign}}\,\lambda|\lambda|^{-z}\dim E_{\lambda}(D^+_{\Ec}),
	$$
	where $E_{\lambda}(D^+_{\Ec})$ is the eigenspace of $D_{\Ec}^+$
	associated to $\lambda$. This extends meromorphically to the whole
	complex plane with the origin not being a pole and $\eta^+_{\Ec}(0)$
	is a well defined real number called the {\it eta invariant} of
	$D_{\Ec}^+$. It measures the overall asymmetry of ${{\rm
			Spec}}(D_{\Ec}^+)$ with respect to the origin.

	In general, the existence of a boundary implies that the space of
	harmonic spinors
	$\ker\partial_{\Ec}=\ker\partial_{\Ec}^+\oplus\ker\partial_{Ec}^-$
	is infinite dimensional and one has to impose suitable boundary
	conditions in order to restore finite dimensionality of kernels.
	Here we consider Atiyah-Patodi-Singer (APS) type boundary conditions
	and for this we need to introduce some notation. If $\Dc$ is a self
	adjoint elliptic operator acting on sections of a bundle
	$\Fc\to\Sigma$, we denote by $\Pi_{I}(\Dc):L^2(\Fc)\to L^2(\Fc)$ the
	spectral projection of $\Dc$ associated to the interval
	$I\subset\mathbb{R}$. Also, if $\psi\in\Gamma(S_W\otimes\Ec)$ we set
	${\varphi}=\psi|_{\Sigma}$. Now  consider $\Gamma_{\geq
		0}(S^+_W\otimes\Ec)=\{\psi\in
	\Gamma(S^+_W\otimes\Ec);\Pi_{[0,+\infty)}(D_{\Ec})\varphi=0\}$
	and $\Gamma_{
		>0}(S^-_W\otimes\Ec)=\{\psi\in \Gamma(S^-_W\otimes\Ec);
	\Pi_{(0,+\infty)}(D_{\Ec})\varphi=0\}$, which are the domains of
	the operators
	\[
		\dirac^+_{\Ec,\geq 0}={{\dirac^+_{\Ec}}|}_{\Gamma_{\geq
				0}(S^+_W\otimes\Ec)}:\Gamma_{\geq 0}(S^+_W\otimes\Ec)\to
		\Gamma(S^-_W\otimes\Ec) \label{dec1}
	\]
	and
	\[
		\dirac^-_{\Ec,> 0}={{\dirac^-_{\Ec}}|}_{\Gamma_{>
				0}(S^-_W\otimes\Ec)}:\Gamma_{>0}(S^-_W\otimes\Ec)\to
		\Gamma(S^+_W\otimes\Ec), \label{dec2}
	\]
	respectively. These are adjoints to each other and  moreover
	$\dirac^+_{\Ec,\geq 0}$ is a Fredholm operator with a well defined
	index
	\[
		{{\rm ind}}\,\dirac^+_{\Ec,\geq 0}=\dim\ker \dirac^+_{\Ec,\geq
			0}-\dim\ker \dirac^-_{\Ec, > 0}.
	\]
	The following formula computes this index (see \cite{atiyah1975spectral} or
	\cite{booss1993bavnbek} for the case where $\Uc$ is an isometric product and
	\cite{gilkey1993index} for the general case):
	\begin{equation}\label{indfor}
		{{\rm ind}}\,\dirac^+_{\Ec,\geq 0}=\int_W[\widehat{A}(TW)\wedge{{\rm
				ch}}(\Ec)]_{2k}+\int_{\Sigma}[\Tc\widehat{A}(TW)\wedge{{\rm
				ch}}(\Ec)]_{2k-1}-\xi_{\Ec}^+(0),
	\end{equation}
	where $\Tc\widehat{A}$ is the transgression of $\widehat{A}$ along $\Sigma$,
	which is polynomial in the curvature, and
	$$
	\xi_{\Ec}^+(0)=\frac{1}{2}\left(\eta_{\Ec}^+(0)+\dim\ker
	D_{\Ec}^+\right).
	$$

	For further reference we present now the integral version of (\ref{dec}). Recalling
	that $\varphi=\psi|_{\Sigma}$,  the identity reads
	\begin{equation}\label{weit:twist}
		\int_X\left(\frac{\kappa}{4}|\psi|^2+\langle
		\Rc^{[\nabla]}\psi,\psi\rangle+|\nabla\psi|^2-|\dirac_{\Ec}\psi|^2\right)
		=\int_{\Sigma}\left(\langle
		D_{\Ec}\varphi,\varphi\rangle-\frac{H}{2}|\varphi|^2\right).
	\end{equation}
%	The proof is the same as in \cite{hijazi2000dirac}, where the untwisted version appears. 

	\section{Gromov's $K$-area}
	\label{kar}
	
	The concept of $K$-area was introduced by Gromov \cite{gromov1996positive} in order
	to quantify previous results on geometric-topological obstructions
	to the existence of metrics with positive scalar curvature. We now
	briefly review this material.
	
	Let $(X^{2k},g)$ be a closed Riemannian manifold (not necessarily
	spin). By pulling
	back any of the half-spin bundles $S^\pm_{\mathbb S^{2k}}\to \mathbb S^{2k}$, whose top Chern classes are non-trivial, under
	a degree one map we have that the set of complex vector bundles
	over $X$ which are {\em homologically non-trivial} (i.e. which have
	at least a nonzero Chern number) is nonempty. Notice that by
	Chern-Weil theory, the Chern numbers, which are topological
	invariants of $\Ec$, can be computed by integrating over $X$ certain
	universal differential forms depending on the curvature tensor
	$R^{\Ec}$ of any compatible connection on $\Ec$. Thus $\Ec$ is {\em
		homologically trivial} (i.e. all Chern numbers vanish) if
	$R^{\Ec}=0$. Here and in what follows we abuse notation and denote
	simply by $R^{\Ec}$ the curvature of {\it any} compatible connection
	on $\Ec$.
	
	We then let $(\Ec,\nabla)$ vary over the set of homologically non-trivial
	hermitian bundles (and compatible connections) over $X$ and define
	the $K$-area of $(X,g)$ by
	\begin{equation}\label{karea}
		\ka(X,g)=\sup\frac{1}{\|R^{\Ec}\|_g},
	\end{equation}
	where
	\begin{equation}\label{norm}
		\|R^{\Ec}\|=\sup_{v\wedge w\neq 0} \frac{\|R^{\Ec}_{v,
				w}\|_g}{\,\,\|v\wedge w\|_g}
	\end{equation}
	and $\|v\wedge w\|^2_g=g(v,v)g(w,w)-g(v,w)^2$. Clearly, the $K$-area
	as defined above is a Riemannian invariant but the fact that it is
	finite or infinite is a topological property of $X$. 
	We note that the use of the operator norm in the definition allows us to conclude that if $(\Ec_1,\nabla_1)$ and $(\Ec_2,\nabla_2)$ then
	\begin{equation}\label{oper:cons}
		\left\|R^{\nabla_1\oplus\nabla_2}\right\|_g=\max\left\{\left\|R^{\nabla_1}\right\|_g,\left\|R^{\nabla_2}\right\|_g\right\},\quad 	\left\|R^{\nabla_1\otimes\nabla_2}\right\|_g=\left\|R^{\nabla_1}\right\|_g+\left\|R^{\nabla_2}\right\|_g. 	
	\end{equation}
	
	This is a rather flexible concept which can be easily adapted to other settings. For instance, if
	 $X\in\bg_{2k}$ we may retain the definition
	(\ref{karea}) but restricting to bundles which are trivial at infinity
	(i.e. in a neighborhood of the point at infinity in the one-point
	compactification of $X$). The allowable connections are required to
	be flat at infinity so that characteristic numbers related to $\Ec$
	are obtained by integrating over $X$ characteristic differential
	forms with compact support.  
	By the pull back construction above, we always have that $\ka(X,g)>0$.
	In any case, with this definition, the fact that the $K$-area is
	finite or infinite is obviously a bounded geometry type invariant of
	$X\in\bg_{2k}$, so we may simply write $\ka(X)=\ka(X,g)$ whenever no confusion arises. 
	It turns out, however, that this notion of $K$-area in $\bg_{2k}$ is not suitable for the purposes we have in mind; see Remark \ref{non:suit}. In any case,  
	the following useful characterizations are readily
	derived from the definitions.
	
	\begin{proposition}\label{grom}
		i) $\ka(X,g)=+\infty$ if and only if for any $\epsilon>0$ there
		exists a homologically non-trivial admissible $\Ec$ over $X$ with
		$\|R^{\Ec}\|_g\leq\epsilon$; ii) $\ka(X,g)<+\infty$ if and only if
		there exists $\epsilon_{X,g}>0$ such that if $\Ec$ over $X$ is
		admissible and $\|R^{\Ec}\|_g\leq\epsilon_{X,g}$ then $\Ec$ is
		homologically trivial.
	\end{proposition}

	\begin{remark}\label{rank:blows}
		It is crucial here to let $\rm rank\,\Ec\to+\infty$ as  $\epsilon\to
		0$ in case ${\rm vol}_{2k}(X)<+\infty$. In fact, by Chern-Weil theory we have $1\leq
		\|R^{\Ec}\|_g^k{{\rm vol}}_{2k}(X)|P({{\rm rank\,\Ec}},k)|$ for
		some polynomial $P$ in case $\Ec$ is homologically non-trivial.
	\end{remark}
	
	\begin{remark}\label{hyper}
		Examples of manifolds with infinite $K$-area in the closed category include tori and, more generally, finitely enlargeable manifolds \cite{gromov1980spin}. This includes solvable manifolds and non-positively curved manifolds whose fundamental group is residually finite. On the other hand, 
		simply connected closed manifolds
		always have finite $K$-area. In the bounded geometry category, let us take a  Riemannian
		manifold $X^{2k}$ which is {\em hyper-euclidean}  in the sense that for
		some $\epsilon>0$ there exists a nonzero degree, proper map
		$f:X\to\mathbb{R}^{2k}$ satisfying $\|f_*\omega\|\leq
		\epsilon\|\omega\|$ for any $2$-form $\omega$ over $X$. This
		gives $\ka(X)\geq\epsilon^{-1}\ka(\mathbb{R}^{2k})$ and
		sending $\epsilon\to 0$ (which can be accomplished by scaling) we
		get $\ka(X)=+\infty$. As a simple example we can
		take $X$ to be any Hadamard manifold  within its $\bg_{2k}$ class. In particular, $\ka(\mathbb{R}^{2k})=+\infty$. 
	\end{remark}

	\section{Uniformly finite homology, bounded de Rham cohomology and $K_{\rm area}^b$}
	\label{bd}
	
	As  mentioned in the Introduction, Block and Weinberger \cite{block1992aperiodic} have
	defined a bound\-ed geometry (in fact, coarse) invariant homology
	$H^{\rm uf}_0$, the so-called {\it uniformly finite homology} in degree
	zero. If $Y\in\bg_n$ then $S\subset Y$ defines a class $[S]\in
	H^{\rm uf}_0(Y)$ if $S$ is locally uniformly finite. 	Our aim here is to extend Whyte's obstruction to $Y\sharp_SM$ in case the attached manifold satisfies $K_{\rm area}(M)=+\infty$. For this we need the appropriate notion of $K$-area for objects in $\bg_{2k}$.

		Now, a key property of $H^{\rm uf}_0$ is that it is
	naturally dual to $H^n_b$, the {\it bounded de Rham cohomology} in
	degree $n=\dim X$ \cite[Lemma 2.2]{whyte2001index}. If $\Upsilon\in H^n_b(X)$ we denote by
	$\Upsilon^{\rm uf}$ the corresponding homology class in $H_0^{\rm uf}(X)$. Regarding this identification, the following is well-known; see
	\cite[Theorem 3.1]{block1992aperiodic}.
	
	\begin{proposition}\label{bg:known:0}
		$X\in \bg_n^N$ if
		and only $X$ is open at infinity, which means by definition that
		large compact domains in $Y$ satisfy a linear isoperimetric inequality. Equivalently, any bounded $n$-form is the differential of a bounded $(n-1)$-form.
	\end{proposition}
	
	The next result goes one step further and provides a useful criterion to decide when a given class in 	$H^n_b(X)\simeq H_0^{\rm uf}(X)$ is trivial; see \cite[Lemma 2.4]{whyte2001index}.

	\begin{proposition}\label{bg:known}		
		$\Upsilon\in H^n_b(X)$  vanishes if and only if for all sufficiently large compact domain $W\subset X$ with boundary $\Sigma$ there holds
		\[
		\left|\int_W\Upsilon\right|\leq C{\rm vol}_{2k-1}(\Sigma), 
		\]
		where  $C>0$ depends on curvature and second fundamental form bounds (and possibly on $\Upsilon$). 
	\end{proposition} 

	It is also proved in \cite[Lemma 2.1]{whyte2001index} that any characteristic polynomial $\mathfrak P$ gives, after evaluation on the curvature, a characteristic number $\mathfrak P^b\in H^n_b(X)$.	
A closer
look at the proof of the Poincar\'e duality $H^{\rm uf}_0\cong H^n_b$ in
\cite{whyte2001index} reveals that this construction behaves quite well under  the infinite
connected sum operations we are dealing with.

\begin{proposition}\label{well}
	If $\Upsilon\in H^n_b(Y)$ is given by a characteristic form via
	Chern-Weil theory, $[S]\in H^{\rm uf}_0(Y)$ and $M$ is closed then
	\begin{equation}\label{welll}
		\Upsilon^{\rm uf}(Y\sharp_S M)=\Upsilon^{\rm uf}(Y)+\Upsilon(M)[S],
	\end{equation}
	where $\Upsilon(M)=\int_{M}\Upsilon$ is the corresponding characteristic
	number computed over $M$.
\end{proposition}

The following result is proved in \cite[Theorem 2.3]{whyte2001index}.

\begin{theorem}\label{why:obs}
	If $X\in \bg_n$ is spin and carries a metric with $\kappa\geq 0$ then $\widehat A^b(X)=0$. 
\end{theorem}

\begin{remark}\label{non:suit}
	This last result is combined with (\ref{welll}) in \cite{whyte2001index} to yield the obstruction mentioned in Remark \ref{exam1}. Thus, 
it is tempting to try to extend Theorem \ref{why:obs} by replacing $\widehat A^b(X)\neq 0$ with  infinite $K$-area as defined in Section \ref{kar}. However, counterexamples to the corresponding statement are easily found. For instance, the bounded geometry class of the standard flat metric in $\mathbb R^{2k}$ carries a scalar flat metric but has infinite $K$-area by Remark \ref{hyper}. 
\end{remark}

The counterexample in the previous remark asks for the proper modification of the notion of $K$-area in $\bg_{2k}$. With this goal in mind, 
we consider Hermitean bundles $(\Ec,\nabla)$ over $X\in\bg_{2k}$ such that:
\begin{itemize} 
	\item the curvature $R^{\nabla}$ is uniformly bounded over $X$, so the corresponding characteristic numbers lie in $H^n_b(X)$.
%	\item the rank $r={\rm rank}\,\Ec$ is uniformly bounded by some sufficiently large constant.
\end{itemize} 
Notice that we do {\em not} require that the bundle $\Ec$ is trivial in a neighborhood of infinity.
Thus, our choice here departs a bit from the general philosophy in \cite{gromov1996positive}, where the allowable bundles in the open case are tied to compactly supported cohomology; see Section \ref{kar}. In any case, with this notion of {\em admissible} bundle at hand, 
the idea is to retain the definition  (\ref{karea}) but restricting to Hermitian bundles $(\Ec,\nabla)$ which are $b$-{\em homologically non-trivial} in the sense that at least one {\em bounded} Chern number $H^{2k}_b(X)\ni c_{I}^b({\Ec})\neq 0$, $|I|=k$. We denote the corresponding geometric invariant by $K^b_{\rm area}(X)$, the  $K^b$-area of $X$.

It might well happen that the set of $b$-homological non-trivial bundles is empty, in which case the bounded $K$-area vanishes for trivial reasons. 
By Proposition \ref{bg:known:0}, this happens if $X$ is non-amenable. Thus, it only possibly holds that $0<K_{\rm area}^b(X)\leq +\infty$ if $X$ is amenable. In any case, the fact that this invariant vanishes, is finite or infinite is a bounded geometry invariant property of $X$. 

The next property helps to detect an important class of manifolds with null $K^b_{\rm area}$
	
	\begin{definition}\label{open:area}
	We say that $X\in\bg_n$ is {\em large at infinity} if there exists an exhaustion of $X$ by  regular compact domains $W$ and $c>0$ such that ${\rm vol}_{n-1}(\Sigma=\partial W)\geq c$ as $W$ sweeps out $X$. 
\end{definition}

\begin{proposition}\label{large:triv}
	If $X\in\bg_{2k}$ is large at infinity and there exists a neighborhood $U$ of infinity with the  homotopy type of the sphere $\mathbb S^{2k-1}$ then $K_{\rm area}^b(X)=0$. In particular, $K_{\rm area}^b(\mathbb R^{2k})=0$. 
\end{proposition}

\begin{proof}
	By Bott periodicity, any (admissible) bundle $\Ec$ over $U$ is stably trivial, that is, 
	there exists $l$ such that
	$\Ec'=\Ec\oplus \Theta^l$ is trivial over ${U}$ 
	(here and in the following,
	$\Theta^l$ is the trivial bundle with typical fiber $\mathbb{C}^l$ and varying base space; this trivial bundle is always assumed to be  endowed with a
	flat connection).
	We note that $l$ may be chosen to depend only on $k$ and $r$. Hence, each Chern number $c_{I}^b({\Ec'})$ is induced by a compactly supported form. Since both the curvature and the rank of $\Ec'$ are controlled, we have  
	\[
	\frac{\left|\int_Wc_{I}(\Ec')\right|}{{\rm vol}_{2k-1}(\Sigma)}\leq \frac{C{\rm vol}_{2k}\left({\rm supp}\,c_I(\Ec')\right)}{c},
	\]
	so  $c_I^b(\Ec')$ vanishes in $H^{2k}_b(X)$ by Proposition \ref{bg:known}. Since $\Ec$ and $\Ec'$ have the same bounded Chern numbers, 
	we see that any such $\Ec$ is $b$-homologically trivial.
\end{proof}

	\section{The proofs of Theorems \ref{aux2} and \ref{aux1}}
	\label{proofaux1}

Proposition \ref{large:triv} is a bit discouraging as it shows that elements in a large class of amen\-able manifolds have null $K^b$-area for trivial reasons (no $b$-homological non-trivial bundle is available).	
Fortunately, starting with {amenable} manifolds with null 
or finite  $K^b$-area, the infinite connected sum construction above allows us to display a large supply of manifolds with infinite $K^b$-area. This is precisely the content of Theorem \ref{aux2}, whose proof we now present.

\subsection{The proof of Theorem \ref{aux2}}\label{aux2:proof}
	Let $X=Y\sharp_S M$ and assume by
	contradiction that $0\leq K_{\rm area}^b(X)<+\infty$. Thus, if $g$ is a metric in the given
	bounded geometry class, by the analogue of Proposition \ref{grom} there exists
	$\epsilon>0$ such that if $(\Ec,\nabla)$ admissible over $X$ satisfies
	$\|R^{\nabla}\|\leq \epsilon$ then $\Ec$ is $b$-homologically
	trivial. 
	
	Let $(\mathcal F,\nabla_{\mathcal F})$ an admissible bundle over $Y$ with $\|R^{\nabla_{\mathcal F}}\|\leq\epsilon'< \epsilon$. If 
	$\mathcal{V}\subset M$ is a compact tubular neighborhood of the
	sphere $\mathbb S^{2k-1}$ over which the connected sum operation leading to
	$X$ was carried out, then $\mathcal{V}$ has the same homotopy type
	as $\mathbb S^{2k-1}$. Hence, by Bott periodicity, there exists $l$ such that
	$\Fc'=\Fc\oplus \Theta^l$ is trivial over $\mathcal{V}$. In particular, $\Fc'$ may be extended to $M$ (as a trivial bundle) and hence  to  $X$. 
	
	On the other hand, since $\ka(M)=+\infty$, there exists a homologically non-trivial
	bundle $(\Gc,\nabla_{\Gc})$ over $M$ with $\|R^{\nabla_\Gc}\|\leq \epsilon'$. 
	Again, there exists $m$ such that $\Gc'=\Gc\oplus\Theta^m$ is trivial if restricted to $\mathcal V$, so $\Gc'$ may be extended to $Y$ (as a trivial bundle) and hence to $X$. As before, $l$ and $m$ may be chosen to depend only on $k$ and $r$.
	
	We next consider the bundle $\Hc=\Fc'\oplus\Gc'$.  Since we assume that $0\leq K_{\rm area}^b(X)<+\infty$, $\Hc$ is $b$-homologically trivial over $X$ and by 	duality we find that $c_I(\Hc)^{\rm uf}(X)=0$ for {\em any} bounded Chern number class
	$c_I^b(\Hc)$. However, we may find $J$ such that $c_J(\Hc)(M)=c_J(\Gc')(M)=c_J(\Gc)(M)\neq 0$. 
	By
	Proposition \ref{well},
	$c_J(\Hc)^{\rm uf}(Y)=-c_J(\Hc)(M)[S]\neq 0$, so $\Hc|_Y$ is
	$b$-homologically non-trivial. Since $\Hc|_Y$ and $\Fc$ have the same bounded Chern numbers, $\Fc$ is $b$-homologically non-trivial as well.  But this means that
	$K_{\rm area}^b(Y)=+\infty$, a contradiction that completes the proof of Theorem \ref{aux2}.

\subsection{The proof of Theorem \ref{aux1}}\label{aux1:proof}
	We fix  $X\in\bg_{2k}$ spin with $K^b_{\rm area}(X)=+\infty$ and consider $W\subset X$ a compact regular domain with $\Sigma=\partial W$. 
	 By passing to another metric in the bounded geometry class if needed, we may assume that $\kappa>0$ on $W$. 
	Also, we may assume that the second fundamental form of $\Sigma$ is uniformly bounded. Note that the existence of arbitrarily large domains satisfying this latter condition with bounds depending on the underlying geometry follows from the chopping construction in \cite{cheeger1990chopping} combined with the bounded geometry assumption. 
	
		We first observe that applying index theory to $W$ and arguing as in the proof of  
\cite[Theorem 2.3]{whyte2001index}, we obtain
	\begin{equation}\label{why:est:unt}
	\left|\int_W[\widehat A(TX)]_{2k}\right|\leq C {\rm vol}_{2k-1}(\Sigma),
	\end{equation}
	which means that $\widehat A^b(X)=0$ by Proposition \ref{bg:known}. Henceforth, $C>0$ is a constant which only depends on the dimension and the curvature bounds (in case it further depends on additional parameters, we explicitly indicate this by using them as subscripts). 
	Our aim is to adapt to the twisted case the argument leading to (\ref{why:est:unt}).
	
	Let $\kappa_W:=\inf_W\kappa_g>0$. By the analogue of Proposition \ref{grom}, for each $\epsilon>0$
	there exists an admissible $b$-homologically non-trivial bundle
	$(\Ec,\nabla)$ over $X$ with $\|R^{\nabla}\|_{\kappa_W g}\leq
	\epsilon$, or equivalently, $\|R^{\nabla}\|_{g}\leq
	\epsilon\kappa_W$. In view of (\ref{jos3}) this gives
	$\|\Rc^{[\nabla]}\|_g\leq \epsilon C_{k}\kappa_W$ for some
	$C_{k}>0$ depending only on $k$. It follows that, restricted to $W$,
	\[
	\Rc^{[\nabla]}+\frac{\kappa}{4}\geq \left(-\epsilon
	C_{k}+\frac{1}{4}\right)\kappa_W,
	\]
	so  if $\epsilon=\epsilon_{k,W}\leq\min \{1/8C_{k},k_W^{-1}\}$ we obtain not only the curvature bound 
	$	\|R^{\nabla}\|_g\leq 1$, but also 
	the pointwise estimate
	\begin{equation}\label{point}
		\Rc^{[\nabla]}+\frac{\kappa}{4}\geq \frac{\kappa_W}{8}>0
	\end{equation}
	on $W$.

	As in the proof of (\ref{why:est:unt}) we next appeal to the index formula (\ref{indfor}). 
	For simplicity we set $\mathsf s={\rm vol}_{2k-1}(\Sigma)$ and $r={\rm rank}\,\Ec$. It is well-known that 
	$|\xi_{\Ec}^+(0)|\leq Cr\mathsf s$ \cite[Theorem 3.1.1]{ramachandran1993vonneumann}. 
	Moreover, from the curvature bounds  at our disposal  we clearly have 
	\begin{equation}\label{est:prel}
		\left|\int_{\Sigma}[\Tc\widehat{A}(TX)\wedge{{\rm
				ch}}(\Ec)]_{2k-1}\right|\leq C\mathsf P(r)\mathsf s,
	\end{equation}
where $\mathsf P$ is a certain polynomial; compare with Remark \ref{rank:blows}.

	We now examine the index term in (\ref{indfor}). We take a harmonic spinor $\psi\in {{\rm ker}}\,\dirac_{\Ec,\geq
		0}^+\cup {{\rm ker}}\,\dirac_{\Ec,>
		0}^-$. By (\ref{point}) and (\ref{weit:twist}), 
	\[
	\int_{\Sigma}\langle D_{\Ec}\varphi,\varphi\rangle\geq a\int_{\Sigma} |\varphi|^2,\quad \varphi=\psi|_{\Sigma}, 
	\]
	where we may assume that $a=\inf_{\Sigma_i} H/2<0$. 
	We now Fourier expand $\varphi=\sum_{\lambda\leq 0}c^{\varphi}_{\lambda}\phi_{\lambda}$, where $D_{\Ec}\phi_{\lambda}=\lambda\phi_{\lambda}$. The left-hand side then becomes $\sum_{\lambda\leq 0}\lambda|c^\varphi_{\lambda}|^2$
	and we see at once that at least one $\lambda$ lies in the interval $[a,0]$ if $\psi\neq 0$. In other words, the linear map 
	\[
	\Xi: {{\rm ker}}\,\dirac_{\Ec,\geq
		0}^+\cup {{\rm ker}}\,\dirac_{\Ec,>
		0}^-\to \Pi_{[a,0]}(D_{\Ec})(L^2(\Ec|_{\Sigma})), \quad 
	\Xi(\psi)= \sum_{\lambda\in[a,0]}c^\varphi_\lambda\phi_\lambda, 
	\]
	is well defined and injective. 
	The conclusion is that 
	\[
	{\rm dim}\,{{\rm ker}}\,\dirac_{\Ec,\geq
		0}^++{\rm dim}\, {{\rm ker}}\,\dirac_{\Ec,>
		0}^-\leq N(a^2),
	\]
	where $N(a^2)$ is the number of eigenvalues of the Dirac Laplacian $D_{\Ec}^2$ in the interval $[0,a^2]$. A standard argument \cite{whyte2001index,bordoni1994spectral} shows that $|N(a^2)|\leq Cr\mathsf s$. As a consequence,  $|{\rm ind}\,\dirac_{\Ec,\geq
		0}^+|\leq Cr\mathsf s$.

		\begin{remark}\label{rank:rem}
		The reason why the estimates on $|\xi^+_{\Ec}(0)|$ and $|{{\rm ind}}\,\dirac^+_{\Ec,\geq 0}|$ are {\rm linear} in $r$ is that the corresponding arguments are spectral in nature and rely on the cuvature bounds to compare the spectra of $D_\Ec$ and $D_{\Theta^r}$, where $\Theta^r=\Sigma\times\mathbb C^r$ is the trivial bundle endowed with a flat connection, and then using that this latter Dirac operator is just a sum of $r$ copies of $D$. On the other hand, since we are not assuming that $\Ec|_{\Sigma}$ is trivial, the estimate in (\ref{est:prel}) is polynomial in principle.
	\end{remark}

	 If we now combine the estimates above  with (\ref{indfor}), the conclusion is that 
	 \begin{equation}\label{est:pr:gen}
	 \left|\int_{W}[\widehat{A}(TX)\wedge{\widehat{\rm ch}}(\Ec)]_{2k}\right|	
	 = \left| \int_{W}[\widehat{A}(TW)\wedge{{\rm ch}}(\Ec)]_{2k}- r\int_{W}[\widehat{A}(TX)]_{2k}  \right|
	 \leq C_r\mathsf s.
	 \end{equation}
	 Here,  the {\em reduced} Chern character of $\Ec$ is 
	 \[
	 {\widehat{\rm ch}}(\Ec)={\rm ch}(\Ec)-r
	 ={\rm ch}_1(\Ec)+{\rm ch}_2(\Ec)+\ldots,
	 \]
	 with ${\rm ch}_i(\Ec)\in \Gamma(\wedge^{2i}TX)$  defined by a universal homogeneous characteristic polynomial of degree $i$   in  $R^{\nabla}$.

		We now follow \cite{gromov1996positive} and bring the mechanism of Adams operations to our discussion. Recall that this is a rule that to each $\mu\in \mathbb N$  and $\Ec$ as above associates a (virtual) bundle $\Psi_\mu\Ec$  
	which is a universal expression in terms of tensor products of exterior powers of $\Ec$. It is compatible with the Chern character map in the sense that 
	\[
	{\rm ch}(\Psi_\mu\Ec)=\sum_{j\geq 0}{\rm ch}_j(\Ec)\mu^j. 
	\]
	In particular, ${\rm rank}\,\Psi_\mu\Ec=r={\rm rank}\,\Ec$.
	%	, so that 
	%		\[
	%	\widehat{\rm ch}(\Psi_\mu\widetilde\Ec)=\sum_{j\geq 1}{\rm ch}_j(\widetilde \Ec)\mu^j.
	%	\]
	Moreover, for each $\nu\in\mathbb N$ and a multi-index $\mu_{(\nu)}=(\mu_1,\ldots,\mu_\nu)$, one has that 
	\[
	\Psi_{\mu_{(\nu)}}\Ec:=\Psi_{\mu_1}\Ec\otimes\cdots\otimes \Psi_{\mu_\nu}\Ec
	\] 
	satisfies ${\rm rank}\,\Psi_{\mu_{(\nu)}}\Ec=r^\nu$ and 
	\[
		{\rm ch}_j\Psi_{\mu_{(\nu)}}\Ec=\sum_{i_1+\ldots+i_\nu=j}\mu_1^{i_1}\cdot\ldots\cdot \mu_\nu^{i_\nu}{\rm ch}_{i_1}(\Ec)\wedge\ldots\wedge {\rm ch}_{i_\nu}(\Ec). 
	\] 
Also, by (\ref{oper:cons}) the curvature bounds are `'stable'' under Adams operations, 
	so
the argument leading to (\ref{est:pr:gen}) works fine 
to yield
\[
 \left|\int_{W}[\widehat{A}(TX)\wedge{\widehat{\rm ch}}(\Psi_{\mu_{(\nu)}}\Ec)]_{2k}\right|\leq C_{r,\nu}\mathsf s.
\]
In particular, if we take $\nu=k$ then Proposition \ref{bg:known} implies that $\mathscr Q^b(\mu_{(k)})=0$, where
	 \begin{eqnarray*}
	 	\mathscr Q(\mu_{(k)})
	 	&=& [\widehat{A}(TX)\wedge{\widehat{\rm ch}}(\Psi_{\mu_{(k)}}\Ec)]_{2k} \\
	 	& = & 
	 	\sum_{i_1+\ldots+i_k=k}\mu_1^{i_1}\cdot\ldots\cdot \mu_k^{i_k}{\rm ch}_{i_1}(\Ec)\wedge\ldots\wedge {\rm ch}_{i_k}(\Ec)+\ldots,
	 \end{eqnarray*}
	 with the dots corresponding to `'lower order terms''. 
	 Thus, each bounded Chern character class $({\rm ch}_{i_1}(\Ec)\wedge\ldots\wedge {\rm ch}_{i_k}(\Ec))^b$ vanishes. Since  each bounded Chern number  ${c}^b_I(\Ec)$ is a universal rational linear combination of such classes, we see that each  ${c}^b_I(\Ec)$ vanishes as well.
But this means that $\Ec$ is $b$-homologically trivial and this contradiction completes the proof of Theorem \ref{aux1}.

	\bibliographystyle{alpha}
	\bibliography{karea-inf-conn-arxiv-v1}

\end{document}